\documentclass{article}
\usepackage{fleqn,times,cite}
\usepackage[intlimits]{amsmath}
\usepackage{amsfonts}
\usepackage{amssymb}
\usepackage{hhline}

\twocolumn

\newtheorem{defn}{Definition}
\newtheorem{theorem}{\sf\bfseries Theorem}
\newtheorem{lemma}{\sf\bfseries Lemma}

\newtheorem{coro}{\sf\bfseries Corollary}
\def\IR{{\rm I\!R}}
\def\IC{{{\rm C}\!\!\!{\sf I}~}}
\def\QED{~\rule[-1pt]{5pt}{5pt}\par}

\newenvironment{mat}{\left[\begin{array}}{\end{array}\right]}

\renewcommand{\t}{^{\mbox{\tiny\sf T}}}
\newcommand{\dt}{{\mbox{\rm\em dt}}}

\newcommand{\hs}{\hspace{4mm}}
\newcommand{\hsm}{\hspace{-1.4mm}}
\newcommand{\he}{{\rm He}}

\newcommand{\bfOme}{{\bf\Omega}}

\newcommand{\Ltwo}{{{\cal L}_2[0,\infty)}}

\renewcommand{\j}{{\jmath}}

\newcommand{\inte}{{\rm Int}}

\newcommand{\mylabel}[1]{\label{#1}}
\newcommand{\mR}{\mathbb R}
\newcommand{\mZ}{\mathbb Z}

\title{\large\bf CONIC S-PROCEDURE AND CONSTRAINED DISSIPATIVITY
\\
\thanks{
 The work is supported in part by the Russian Foundation Foundation for Basic Research
(project 05-01-00869) and CNRS-RAS cooperative research program
(project 04-16394)} }

\author{~~~ Alexander L. Fradkov\thanks{
A. L. Fradkov is with the Institute for Problems of Mechanical
Engineering, Control of Complex Systems Laboratory, V.O., Bolshoy, 61,
St. Petersburg, 199178, Russia, {\rm alf@control.ipme.ru}}
}

\begin{document}
\date{}
\maketitle

\begin{abstract}
A new version of classical S-procedure in system theory is
proposed based on duality in the space of positive definite
matrices and introduction of matrix Lagrange multipliers. A new
proof and extension of the recent results \cite{IHF05}
%of T.Iwasaki, S. Hara, A. Fradkov ( Systems & Control Letters, 2005. Vol 54 (7), pp 681-691)
concerning equivalence between frequency domain inequality on
finite frequency range and constrained dissipativity property for
linear systems is given. The results of this paper  extend
S-procedure to allow for analysis and design of robust systems
with matrix inequalities constraints.
\end{abstract}

{\bf Keywords.} Kalman-Yakubovich-Popov (KYP) lemma, S-procedure,
frequency domain inequality, linear matrix inequality,
dissipativity.

\section{Introduction} \mylabel{sec:intro}

Recently a number of new tools for systems analysis and design
 related to frequency domain inequalities (FDI) over a finite frequency range
(so called {\it Generalized KYP-lemma}) have been developed
\cite{iwasaki:kyp,iwasaki:dyn,iwasaki:gkyp}. It follows from the
results of \cite{iwasaki:kyp,iwasaki:dyn,iwasaki:gkyp} that
fulfillness of a standard FDI in a finite frequency range is
equivalent to validity of some nonclassical linear matrix
inequalities (LMI) for a pair of matrices $P,Q$ replacing
inequalities for a single matrix $P$ appearing in the classical
KYP-lemma. It was shown in \cite{IHF05} that FDI, in turn, are
equivalent to some time-domain inequality (TDI, dissipation
inequality \cite{willems:72II}), valid only over a part of the
system trajectories, determined by an additional integral matrix
inequality (restricted or constrained dissipativity \cite{IHF05}).
Thus, a complete extension of the classical KYP-results on
equivalence between FDI, TDI and LMI to the ''finite-frequency"
case was obtained. Note that the proof of equivalence between FDI
and TDI in \cite{IHF05} goes along the
 lines of the necessity proof for the frequency-domain absolute stability criterion
\cite{yak73,yakfr73}.

In this paper a new proof of the result of \cite{IHF05} is
provided based on the losslessness result for a new version of the
classical S-procedure \cite{yakubovich:71}. A new version of the
S-procedure, also included in the paper, deals with constraints in
LMI form, or more generally, conic inequalities in linear spaces.

In the next section a new S-procedure results are presented. In Section 3 they are
 applied to the proof of equivalence between TDI and LMI.

We use the following notation. The set of square integrable
functions on $[0,\infty)$ is denoted by $\Ltwo $, $M^\dagger$,
where $M$ is a matrix, stands for  its transposition and complex
conjugate of all elements. For a square matrix $M$, its Hermitian
part is defined by $\he(M):=(M+M^\dagger)/2$.
%For a Hermitian matrix $M$, its
%maximum eigenvalue is denoted by $\lambda_{\max}(M)$. The
%operators $\Re(\cdot)$ and $\Im(\cdot)$ take the real and
%imaginary parts of the arguments, respectively. The set of
%positive integers up to $n$ is denoted by $\I_n$.
 The interior of a set $\bfOme$ is denoted by $\inte~\bfOme$.

\section{Conic S-procedure}
Let $X$, $Y_1$, $\dots, Y_m$ be linear topological spaces, $G_j\!:X\to\! Y_j$,
$j=1,\dots,m$ be continuous mappings.

Let for any $j=1,\dots,m$ a convex cone $K_j\subset Y_j$ be given
defining inequality $G_j(x)\ge 0$ for $x\in X$ as inclusion
$G_j(x)\in K_j$. Let $Y^*_j$ denote an dual space to $Y_j$, i.~e.
a linear space of linear continuous functionals $y^*_j$ on $Y_j$
and $K^*_j\!\subset\! Y^*_j$ denote an dual cone to $K_j$, i.~e.
$K^*_j\!=\!\big\{y_j^*\in Y_j^*: \big<y_j^*,y_j\big>\ge 0~ \forall
y_j\in K_j\big\}$, where $\big<y_j^*,y_j\big>$ is the value of the
functional $y_j^*$ at the element $y_j$.

Obviously, if $Y=\mR^1\times Y_1\times \dots \times Y_m$, then
$Y^*= \mR^1\times Y^*_1\times\dots\times Y^*_m$ is the set of all
corteges $\big(y_0^*,y_1^*,\dots,y_m^*\big)$, where $y_0^*\in
\mR^1$, $y_j^*$ is a linear functional from $Y_j^*$.

Consider the following two relations for the mappings $F_0$,
$G_1$, \dots, $G_m$.
\begin{itemize}
\item[({\it A})] $F(x)\ge 0$ for $x\in X$, $G_j(x)\in K_j$,
$j=1,\dots,m$; \item[({\it B})]$\exists\tau_0\ge 0$, $\tau_j\in
K^*_j$:
$\tau_0F(x)-\sum\limits_{j=1}^{m}{\big<\tau_j,G_j(x)\big>\ge 0} ~
 \forall x\in X$.
\end{itemize}
Obviously, validity of ({\it B}) with $\tau_0 > 0$ implies ({\it
A}). Indeed, if $x\in X$ satisfies inequalities $G_j(x)\in K_j$,
$j=1,\dots,m$, then it follows from ({\it B}) that $\tau_0F(x)\ge
0$, since $\big<\tau_i,G_i(x)\ge 0\big>$ for $j=1,\dots,m$. The
opposite statement is not true even in the case of scalar
constraints $Y_j=\mR^1$, $j=1,\dots,m$, corresponding to the
classical $S$-procedure \cite{yakubovich:71}.

Similarly to the classical case we will say that $S$-procedure
with conic constraints $G_j(x)\ge 0$ is lossless, if ({\it B})
with $\tau_0 > 0$ implies ({\it A}).

It is well known \cite{F73} that losslessness of the classical
$S$-procedure is equivalent to the duality theorem in the
corresponding optimization problem. However, the problem is, in
general, nonconvex and only a few classes of functionals $F$,
$G_1$, \dots, $G_m$ are known to possess the losslessness
property.

For example, classical $S$-procedure is lossless, if $m=1$ and
$F$, $G_1$ are quadratic forms on real or complex linear space
$X$. It is also lossless, if $m=2$ and $F$, $G_1$, $G_2$ are
quadratic (Hermitian) forms on the complex linear space $X$.
However, classical $S$-procedure for quadratic forms is, in
general, lossy for $m\ge 2$ in real case and for $m\ge 3$ in
complex case \cite{F73}. A.~Megretski and S.~Treil proved in 1990
\cite{megretski:93} that the classical $S$-procedure is lossless
for all $m\ge 1$, if $F$, $G_1$, \dots, $G_m$ are integral
quadratic forms on ${\cal L}_2(0, \infty)$. V.~Yakubovich extended
this result to a more broad class of quadratic functionals,
forming the so-called {\it S-system} \cite{Yak92}.

Below an extension of the results of \cite{Yak92} to the case of
the $S$-procedure with conic constraints is formulated. Note that
the general formulation of the $S$-procedure with conic
constraints was presented, e.g. in \cite{MatveevYak93,Matveev98},

\begin{theorem}\label{theorem:1}.
Let $n_0=1$ and $\bar K$ is  the closure of the cone $K$ generated
by the set

${\cal F} (X)=\{(F_0(x), F_1(x),\dots, F_m(x)):~x\in X\}. $

If the cone $\bar K$ is convex, then the S-procedure with conic
constraints is  lossless.

If, in addition, constraints $G_j(x)\in K_j$ are regular, namely
$\exists x_0:~G_j(x_0)\in \inte~K_j$, then one can choose
$\tau_0=1$ in (B).
 \end{theorem}

{\it Proof}. Condition (A) implies that $F(x)\ge 0$ for $G_j(x)\in
\inte K_j$, i.e. intersection of the set ${\cal F}(X)$ and the
open cone $D=\{(-y_0, y_1, \dots, y_m): y_0>0, y_j\in \inte K_j,
j=1,\dots,m\}$ is empty: $D\bigcap {\cal F}(X)=\phi$. Therefore,
$D\bigcap\bar K=\phi$. Applying separation theorem for cones, we
obtain that there exists vector
$\tau^*=(\tau_0^*,\tau_1^*,\dots,\tau_m^*)\in Y^*$ such that
$<\tau_0^*,F(x)>+\sum_{j=1}^m<\tau_j^*,G(x)>\ge 0$ for all $x\in
X$ and $<\tau^*,y><0$ for all $y\in D$, i.e.
$<\tau_0^*,y_0>+\sum_{j=1}^m<\tau_j^*,y_j><0$. For any $j=1,\dots
m$ pick up $y_j\in \inte K_j$ and choose sequences
$y_{0_k}\rightarrow 0,~y_{s_k}\rightarrow 0$ as
$k\rightarrow\infty$, such that $y_{0_k}>0, y_{s_k}\in K_s, s\neq
j$. If $k\rightarrow\infty$, then we obtain $<\tau_j^*,y_j>\le 0$,
i.e. $-\tau_j^*\in K_j^*$. The first part of the theorem is
proved.

Taking $x=x_0$ from regularity condition and $y_0\neq 0$ yields
$\tau_0^*,y_0> 0$, i.e. $\tau_0^*>0$. Dividing the inequality (B)
by $\tau_0$, we arrive at the second statement of the theorem. End
of the proof.

 Our next step is to extend the definition of
$S$-system \cite{Yak92} to the case of conic constraints.

\begin{defn}. Let $F_j$, $j=0$, $1$, \dots, $m$ be mappings from
a Hilbert space $\mZ$ to spaces of self-adjoint operators over
corresponding Euclidean space $\mR^{n_j}$, such that $F_j: \mZ\to
SR(n_j\times n_j)$.

We say that $F_0$, $F_1$, \dots, $F_m$ form a S-system if there
exists a subspace $\mZ_0$ and a sequence of linear bounded
operators $T_k: \mZ\to\mZ$, $k=1$, $2$, \dots ~such that
\begin{itemize}
\item[({\it i})] $<T_kz_1,z_2>\to 0$ as $k\to\infty$ for all
$z_1$, $z_2\in \mZ$; \item[({\it ii})] $\mZ_0$ is invariant for
$T_k$ for all $k=1$, $2$, \dots; \item[({\it iii})] $F_j(T_kz)\to
F_j(z)$ as $k\to\infty$ for all $j=0$, $1$, \dots m, $z\in\mZ_0$.
\end{itemize}
\end{defn}

\begin{lemma} Let $F_j$, $j=0$, $1$, \dots, $m$ form $S$-system.
Define the map ${\cal F}:
\mZ\to\prod\limits_{j=0}^{n}{\mR^{n_j\times n_j}}$ by means of the
relation

\begin{align*}{\cal F}(z)&=F_0(z)\otimes F_1(z)\otimes
\dots \otimes F_m(z)\\&\in \mR^{n_0\times
n_0}\otimes\mR^{n_1\times n_1}\otimes\dots\otimes \mR^{n_m\times
n_m}.\end{align*}

Then the closure of the image ${\cal F}(\mZ)$ is a convex set in
$\mR^{n_0^2+n_1^2+\dots +n_m^2}$.\label{lemma:1}
\end{lemma}

In the special case $n_j=1$, $j=1, \dots, m$, Lemma \ref{lemma:1}
coincides with Lemma 1 of the paper \cite{Yak92} and it is proved
similarly to the Lemma 1 of \cite{Yak92}.

{\bf Example 1.} An important series of examples for $S$-systems
is provided by finite family of integral quadratic operators on
the Hilbert space ${\cal L}_2[0,\infty)$ of square integrable
functions with values $z(t)\in \mR^{n_j}$. The mappings are
defined for any $z\in {\cal L}_2[0,\infty)$ as follows:
\begin{equation}
\label{f1} F_j(z)=\he~\int\limits_0^\infty{F_j^{'}
z(t)z^{\dagger}(t)F_j^{''\dagger}\,dt,}
\end{equation}
where $F_j^{'},~F_j^{''}$ are $n_j\times n_j^{'}$ symmetric
matrices.
 In this case the family of the operators $T_k$ can be chosen
as time shifts: $T_k(z)(t)=z(t+k)$, while the subspace $\mZ_0$ can
be chosen as the set of functions with zero initial conditions:
\begin{align*}
\mZ_0=\big\{z(\cdot): z(\cdot)\in {\cal L}_2(0,\infty),
z_2(0)=0\big\}.
\end{align*}

The proof of the S-system property for Example 1 is again similar
to \cite{Yak92}. Note that the cone of positive semidefinite
matrices is selfdual. Therefore S-procedure with conic constraints
determined by functions (\ref{f1}) deals with positive
semidefinite matrix Lagrange multipliers.

Properties of the S-procedure in general case are given by the
following theorem.

\begin{theorem}\label{theorem:2}. $S$-procedure with conic
constraints is lossless for any family of self-adjoint operators
$F_0$, $F_1$, \dots, $F_m$ forming an $S$-system. \end{theorem}

Proof follows immediately from Theorem \ref{theorem:1} and Lemma
\ref{lemma:1}. The result can be extended to the case of equality
constraints and to the case of the  so called {\it generalized
$S$-procedure} introduced in \cite{iwasaki:kyp}.

%%%%%%%%%%%%%%%%%%%%%%%%%%%%%%%%%%%%%%%
\section{Constrained dissipativity}

In this section,  we first present a special case of the
generalized KYP lemma \cite{iwasaki:gkyp},  characterizing FDIs in
the continuous-time setting. Let complex matrices $A$, $B$, $\Pi$,
and real scalars $\varpi_1$, $\varpi_2$ be given.
%Let $\tau$ be $+1$ or $-1$, and
Define
\begin{equation} \mylabel{Ome}
\bfOme:=\{~\omega\in\IR~|~(\omega-\varpi_1)(\omega-\varpi_2)\leq0~\}.
\end{equation} (We may assume $\omega_2>0$ without loss of
generality).

\begin{theorem} \mylabel{thm:ckyp} \cite{iwasaki:gkyp}.
Suppose $\Pi$ is Hermitian matrix, pair $(A,B)$ is controllable,
and $\bfOme$ has a nonempty interior. Then the following
statements are equivalent.
\begin{itemize}
\item[(i)] The frequency domain inequality
\begin{equation} \mylabel{fdic}
\begin{mat}{c} (j\omega I-A)^{-1}B \\ I \end{mat}^*\Pi
\begin{mat}{c} (j\omega I-A)^{-1}B \\ I \end{mat}\leq0
\end{equation}
holds for all $\omega\in\bfOme$ such that $\det(j\omega I-A)\neq0$.
\item[(ii)] There exist Hermitian matrices $P$ and $Q$ such that \\
$Q\geq0$ and the linear matrix inequality
\begin{equation} \label{lmic}
\begin{mat}{cc} \hsm A & \hsm B \hsm \\ \hsm I & \hsm 0 \hsm \end{mat}^*
\begin{mat}{cc} -Q & \hsm P+j\varpi_oQ \hsm \\
\hsm P-j\varpi_oQ & \hsm -\varpi_1\varpi_2Q \hsm \end{mat}
\begin{mat}{cc} \hsm A & \hsm B \hsm \\ \hsm I & \hsm 0 \hsm \end{mat}
+\Pi\leq0
\end{equation}
holds, where $\varpi_o:=(\varpi_1+\varpi_2)/2$.
\end{itemize}
\end{theorem}

%We mention that similar results have also been obtained in
%\cite{iwasaki:kyp,iwasaki:dyn}.
Choosing the parameters
$\varpi_1=\varpi_2=0$ and $\tau=-1$, the set $\bfOme$ becomes the
entire real numbers, and thus statement (i) becomes the FDI for
all frequencies. In this case, the term associated with $Q$ in the
LMI (\ref{lmic}) becomes positive semidefinite, and hence the best
choice of $Q$ for satisfaction of (\ref{lmic}) is $Q=0$. The resulting
LMI with variable $P$ is exactly the same as the one in
%\cite{rantzer:kyp}, and thus Theorem~\ref{thm:ckyp} reduces to
the standard KYP lemma.

The following result extends the result of \cite{IHF05}. It
provides an equivalence between FDI and time domain dissipation
inequality over a restricted class of input signals.

\begin{theorem} \mylabel{thm:tdic}.
Let complex matrices $A$, $B$, $\Pi$, and real scalars $\varpi_1$,
$\varpi_2$ be given and $\bfOme$ be defined by (\ref{Ome}).
Consider the system
\begin{equation} \mylabel{sysc}
\dot x(t)=Ax(t)+Bu(t), \hs t\in[0,\infty),
%\hs x(0)=0
\end{equation}
where $x(t)\in\IC^n$ is the state and $u(t)\in\IC^m$ is the input.
Assume that $(A,B)$ is controllable, $\Pi$ is Hermitian, and
$\bfOme$ has a nonempty interior. Then the following statements
are equivalent.
\begin{itemize}
\item[(i)] The frequency domain inequality (\ref{fdic})
holds for $\omega\in\bfOme$.
\item[(ii)] The time domain inequality
\begin{equation} \mylabel{tdic}
\int_0^\infty
\begin{mat}{c} x \\ u \end{mat}^*\Pi
\begin{mat}{c} x \\ u \end{mat}\dt\leq0
\end{equation}
holds for all solutions of (\ref{sysc}) with $u\in\Ltwo$ such that
$x(0)=0,~x\in\Ltwo$ and
\begin{equation} \mylabel{iqcc}
\he~\int_0^\infty (\varpi_1x+\j\dot x)(\varpi_2x+\j\dot
x)^*\dt\leq0.
\end{equation}
\end{itemize}
\end{theorem}
\medskip

Note that the corresponding result of \cite{IHF05} was obtained
under additional condition of asymptotic stability for the system
(\ref{sysc}) which is not required in the current statement.

 In Theorem~\ref{thm:tdic}, a general
frequency interval $\bfOme$ is considered for the FDI, and this
has translated to the input constraint described by (\ref{iqcc}).
Though the physical meaning of this constraint may be not clear in
 general, it becomes clear for the following  special case.

\begin{coro} \mylabel{coro:tdic}
Let real matrices $A$, $B$, $\Pi$, and a positive scalar $\varpi$
be given. Suppose $\Pi$ is symmetric and consider the system
(\ref{sysc}) where $x(t)\in\IR^n$ is the state and $u(t)\in\IR^m$
is the input. Assume that $(A,B)$ is controllable. Then the
following statements are equivalent.
\begin{itemize}
\item[(i)] The frequency domain inequality (\ref{fdic}) holds for
all $\omega$ such that $|\omega|\leq\varpi$. \item[(ii)] The time
domain inequality (\ref{tdic}) holds for all $u\in\Ltwo$ such that
$x\in\Ltwo$ and
\begin{equation} \mylabel{iqccl}
\int_0^\infty \dot x\dot x\t\dt\leq
\varpi^2
\int_0^\infty x x\t\dt.
\end{equation}
\end{itemize}
Moreover, the above two statements are equivalent when the two
inequalities ``$\leq$'' are replaced by ``$\geq$.''
\end{coro}

Loosely speaking, the first part of Corollary~\ref{coro:tdic}
states that the FDI in the low finite frequency range means that
the system possesses the property (\ref{tdic}) for the input
signals $u$ that  drive the states not too fast (slowly). The
bound on the ``slowness'' is given by $\varpi$ in the sense of
(\ref{iqccl}). The second part of Corollary~\ref{coro:tdic} makes
a similar statement for the FDI in the high frequency range. To
derive Corollary~\ref{coro:tdic} from Theorem~\ref{thm:tdic} one
needs to put $\omega_1=-\omega_2=\varpi$.

%It should be noted that condition (\ref{iqccl}) is a matrix-valued
%integral quadratic constraint, and is coordinate free, i.e. if
%(\ref{iqccl}) holds for some minimal realization of the
%system with states $x$, then it holds for all other minimal realizations
%with states $Tx$ where $T$ is an arbitrary nonsingular matrix.

{\it Proof of Theorem  ~\ref{thm:tdic}.}  In view of Theorem
\ref{thm:ckyp}  it is sufficient to prove equivalence of (ii) to
the condition (ii) of Theorem ~\ref{thm:tdic} (solvability of
matrix inequality  (\ref{lmic})). The result follows from Theorem
\ref{theorem:2} of the previous section with $m=2$ (one matrix
constraint). Denote
%\begin{array}
\begin{align} \mylabel{FzA}
z=&~\begin{mat}{c} x \\ u \end{mat}, ~F(z)=-\int_0^\infty
\begin{mat}{c} x \\ u \end{mat}^\dagger\Pi
\begin{mat}{c} x \\ u \end{mat}\dt, \nonumber
\\
G_1(z)&=-\he~\int_0^\infty (\varpi_1x+\j\dot x)(\varpi_2x+\j\dot
x)^\dagger\dt.
\end{align}
%\end{array}

Obviously, TDI (\ref{tdic}),(\ref{iqcc}) correspond to the
statement (A)  of the conic S-procedure. At the same time the
statement (B) means existence of a $n\times n$-matrix $\tau^*$
from the dual cone $K^*$ to the cone of positive semidefinite
matrices satisfying inequality $F(z)-<\tau^*, G_1(z)>\ge 0
~~\forall z$ or
\begin{equation} \mylabel{FzB}
%z=~\begin{mat}{c} x \\ u \end{mat},
\he\int_0^\infty \left(-
\begin{mat}{c} x \\ u \end{mat}^\dagger\Pi
\begin{mat}{c} x \\ u \end{mat}-(\varpi_1x+\j\dot x)^*\tau(\varpi_2x+\j\dot
x)\right)\dt\ge 0.
\end{equation}
Replacement of $\tau$ by $Q$ and substitution of $\dot x$ from
(\ref{sysc}) transforms (\ref{FzB}) into LMI (\ref{lmic}). To
verify regularity condition of Theorem \ref{theorem:2} take
$x_0=(\bar x,\bar u)$, where $$\bar x(t)=-(A+\mu I)^{-1}B\exp(-\mu
t),~\bar u(t)= \exp(-\mu t), ~\mu>0$$ and $\mu\in\inte \bfOme$.
 Application of Theorem  \ref{theorem:2}  ends the proof.

{\it Remark}. Similar results hold for discrete-time case.

\section{Conclusions}

The property of the system defined by the item (ii) of  Theorem
~\ref{thm:tdic} and the Corollary can be called {\it constrained
dissipativity} or {\it restricted dissipativity}. It is weaker
than standard passivity or dissipativity conditions and may better
reflect specifications for real systems. At the same time the
property of ``slowness" described by the inequality (\ref{iqccl})
leaves enough flexibility to be useful for robustness analysis of
systems.

 The results of the
paper shed new light on the intimate interrelations between
S-procedure and KYP-lemma. They allow to extend classical
S-procedure tool to allow for analysis and design of robust
systems with matrix inequalities constraints.

\end{document}